\documentclass{amsproc}

\usepackage{setspace}
\usepackage[autostyle]{csquotes}
\usepackage{amssymb}
\usepackage{amscd,amsmath}

\newcounter{dummy} \numberwithin{dummy}{section}
 \newtheorem{thm}[dummy]{Theorem}
 \newtheorem{prop}[dummy]{Proposition}
 \newtheorem{lem}[dummy]{Lemma}
 
 \newtheorem{conj}[dummy]{Conjecture}
 
 \newtheorem{dfn}[dummy]{Definition}

 \newtheorem{rem}[dummy]{Remark}

\renewcommand{\leq}{\leqslant}
\renewcommand{\geq}{\geqslant}

\title[Tracial state space with non-compact extreme boundary]{Tracial state space with non-compact extreme boundary}

\subjclass[]{}

\author[Wei Zhang]{\bfseries Wei Zhang}

\address{
Department of Mathematics \\ 
Purdue University   \\ 
150 N. University St., West Lafayette, IN 47907-2067\\
USA}
\email{zhang406@purdue.edu}

\begin{document}

\vspace{18mm}
\setcounter{page}{1}
\thispagestyle{empty}

\begin{abstract}
Let $A$ be a unital simple separable C*-algebra. If $A$ is nuclear and infinite-dimensional, it is known that strict comparison of positive elements is equivalent to $\mathcal{Z}$-stability if the extreme boundary of its tracial state space is compact and of finite covering dimension. Here we provide the first proof of this result in the case of certain non-compact extreme boundaries. Additionally, if $A$ has strict comparison of positive elements, it is known that the Cuntz semigroup of this C*-algebra is recovered functorially from the Murray-von Neumann semigroup and the tracial state state space whenever the extreme boundary of the tracial state space is compact and of finite covering dimension. We extend this result to the case of a countable extreme boundary with finitely many cluster points. 
\end{abstract}

\maketitle

\section{Introduction}
The set of traces of a C*-algebra is a very important invariant of the algebra. For example, in \cite{30} the tracial state space is part of the invariant used to classify unital simple AH-algebras of slow dimension growth. For a C*-algebra $A$, let $T(A)$ be its tracial state space (i.e. the set of normalized finite traces of $A$) and $\partial_e T(A)$ be the extreme boundary of $T(A)$. It is known that $T(A)$ is a Choquet simplex if $A$ is unital (\cite{15}). If $A$ is separable, then $T(A)$ is metrizable. 

In this paper, we will only consider unital simple separable C*-algebras. The tracial state space of such C*-algebras can still be very complicated, such as the Poulsen simplex, in which the extreme points are dense(\cite{17}). In fact, \cite{16} shows that every metrizable Choquet simplex occurs as the tracial state space of some simple unital AF-algebra. 

Several recent results in C*-algebras theory have been obtained under the assumption of a compact extreme boundary of the tracial state space. In 2008 A. Toms and W. Winter made the following conjecture: 
\begin{conj}
Let $A$ be a simple unital nuclear separable C*-algebra. The following are equivalent:\\
(1) $A$ has finite nuclear dimension;\\
(2) $A$ is $\mathcal{Z}$-stable;\\
(3) $A$ has strict comparison of positive elements.
\end{conj}

In 2004, M. R$\o$rdam showed that $\mathcal{Z}$-stability implies strict comparison for unital simple exact C*-algebras(\cite{19}). In 2010, W. Winter proved that finite nuclear dimension implies $\mathcal{Z}$-stability for unital separable simple infinite-dimensional C*-algebras(\cite{24}). H. Matui and Y. Sato proved (3) implies (2) for algebras with finitely many extremal tracial states(\cite{26}). A. Toms, S. White, W. Winter(\cite{14}), E. Kirchberg,  M. R$\o$rdam(\cite{20}) and Y. Sato(\cite{21}) established this result in the case of algebras when the extreme boundary of its tracial state space is non-empty, compact and of finite covering dimension. In the second part of this paper, we will prove the following theorem.

\begin{thm}\label{4a}
Let $A$ be a simple nuclear separable unital infinite-dimensional C*-algebra with non-empty tracial state space. Suppose $\partial_e T(A)=X$ has the tightness property and has finite covering dimension. The following conditions are equivalent:\\
(1) $A$ is $\mathcal{Z}$-stable;\\
(2) $A$ has strict comparison.
\end{thm}
The tightness property in this theorem, to be introduced in the next section, will yield the outcome in the case of certain non-compact extreme boundaries, which has not been reported in previous literature. The methods also allow us to show that the conclusion of Theorem 1.1 holds if, instead of asking for the tightness property, we ask that $A$ has the same tracial state space as a recursive subhomogeneous C*-algebra with finite topological dimension (Proposition \ref{4x}).

M. Dadarlat and A. Toms showed in \cite{1} that for a unital simple separable C*-algebra $A$ with strict comparison of positive elements, the Cuntz semigroup of $A$ is recovered functorially from the Murray-von Neumann semigroup and the tracial state state space $T(A)$ whenever $\partial_e T(A)$ is compact and of finite covering dimension. Although their result can be obtained from the $\mathcal{Z}$-stability result which was mentioned above if $A$ is nuclear, it applies to the nonnuclear case as well. In this paper, we will show that this result still holds if $\partial_e T(A)$ is countable and has only finitely many cluster points (Theorem 5.1).

\section{Preliminaries and Notations}
\subsection{The Cuntz Semigroup}Let $A$ be a C*-algebra and let $\mathcal{K}$ denote the algebra of compact operators on a separable infinite-dimensional Hilbert space. Let \mbox{$(A\otimes \mathcal{K})_+$} denote the set of positive elements in $A\otimes \mathcal{K}$. Given $a,b\in(A\otimes \mathcal{K})_+$, we say that $a$ is Cuntz subequivalent to $b$ (denoted $a\precsim b$) if there is a sequence $(x_n)$ in $A\otimes \mathcal{K}$ such that
\[
\|x_nbx_n^*-a\|\rightarrow 0.
\]
We say that $a$ and $b$ are Cuntz equivalent (denoted $a\sim b$) if $a\precsim b$ and $b\precsim a$. The relation $\precsim$ is clearly transitive and reflexive and $\sim$ is an equivalence relation. 

We define the Cuntz semigroup to be $Cu(A)=(A\otimes \mathcal{K})_+/\sim$, and write $\langle a\rangle$ for the equivalence class of $a\in(A\otimes \mathcal{K})_+$. $Cu(A)$ is indeed an ordered Abelian semigroup when equipped with the partial order
\[
\langle a\rangle\leq\langle b\rangle\Leftrightarrow a\precsim b
\]
and the addition operation
\[
\langle a\rangle+\langle b\rangle=\langle a\oplus b\rangle
\]
using an isomorphism between $M_2(\mathcal{K})$ and $\mathcal{K}$. 
\subsection{Rank Functions}
We denote by $T(A)$ the tracial state space of $A$. Given $\tau$ in $T(A)$, we define a map $d_\tau:A_+\rightarrow[0,+\infty)$ by
\[
d_\tau(a)=\lim_{n\to\infty}\tau(a^{1/n}).
\]
This map is lower semicontinuous, affine and nonnegative. It can be extended natrually to $(A\otimes \mathcal{K})_+$ and we always regard this as its domain. It depends only on the Cuntz equivalence class of $a\in(A\otimes\mathcal{K})_+$. Moreover, it has the following properties:\\
(1)if $a\precsim b$, then $d_\tau(a)\leq d_\tau(b)$;\\
(2)if $a$ and $b$ are mutually orthogonal, then $d_\tau(a+b)=d_\tau(a)+d_\tau(b)$;\\
(3)$d_\tau((a-\varepsilon)_+)\nearrow d_\tau(a)$ as $\varepsilon \to 0$.\\
We then define the rank function of $a\in (A\otimes \mathcal{K})_+$, a map $\iota(a)$ from the tracial state space $T(A)$ to $\mathbb{R}^+\cup\{\infty\}$ given by the formula \mbox{$\iota(a)(\tau)=d_\tau(a)$}. 

\subsection{Strict Comparison and $\mathcal{Z}$-stability}
Let $A$ be a unital C*-algebra. We say that $A$ has strict comparison of positive elements if $a\precsim b$ for $a, b\in (A\otimes \mathcal{K})_+$ whenever
\[
d_\tau(a)< d_\tau(b), \forall\tau\in\{\gamma\in T(A)|d_\gamma(b)<\infty\}.
\]
We say $A$ is $\mathcal{Z}$-stable if $A\otimes \mathcal{Z}\cong A$ where $\mathcal{Z}$ is the Jiang-Su algebra. The Jiang-Su algebra was originally defined in \cite{50} and was constructed as an inductive limit of a sequence of C*-algebras with specified connecting mappings. There are a number of other characterizations of the Jiang-Su algebra. For example, one may view $\mathcal{Z}$ as being the stably finite analogue of the Cuntz algebra $\mathcal{O}_\infty$(\cite{51}, \cite{52}).

\subsection{Choquet Simplices}
We already know that the tracial state space of a unital separable C*-algebra is a metrizable Choquet simplex. For a general metrizable Choquet simplex $K$, given any point $\tau\in K$, there exists a unique Borel probability measure $\mu_\tau$ defined on the extreme boundary $\partial_e K$ such that 
\[
f(\tau)=\int_{\partial_e K} f d\mu_\tau
\]
for any affine continuous function $f$ on $K$ (\cite{18}). We say that $\tau$ is represented by $\mu_\tau$. This result can be extended to affine funtions of first Baire class on $K$, which includes all affine lower semi-continuous functions.(\cite{28}) \\
Moreover, since $K$ is a metrizable simplex, we denote $dist$ a metric on $K$. Hence given two points $\tau$ and $\gamma$ in $K$, the distance between them is $dist(\tau, \gamma)$. Throughout this paper, we will use this notation for the metric on the tracial state space of a unital separable C*-algebra.

\begin{dfn}
Recall that in probability theory, a set $\Gamma$ of Borel probability measures on X is called tight if for every $\varepsilon>0$ there exists a compact subset $F$ of $X$ such that
\[
\mu(F)\geq 1-\varepsilon
\] 
for all $\mu\in\Gamma$.\\
For a unital separable C*-algebra $A$, denote $X$ the extreme boundary of its tracial state space and let $\partial X=\overline {X} \backslash X$. From previous discussion we know that there is a set of Borel probability measures on $X$ representing each $\tau\in\partial X$: $\Gamma_{\partial X}=\{\mu_\tau: \tau\in \partial X\}$. We say that $X$ has the tightness property if $\Gamma_{\partial X}$ is tight.
\end{dfn}
Obviously, if $X$ has only finitely many cluster points, then it has the tightness property.
\subsection{Other Notations}
For convenience, we denote $\mathit{Aff}(K)$ the set of real-valued continuous affine functions on a compact metrizable Choquet simplex $K$, and denote $\mathit{LAff}(K)$ the set of bounded strictly positive, lower-semicontinuous affine functions on $K$, and let $\mathit{SAff}(K)$ be the set of extended real-valued functions which can be obtained as the pointwise supremum of an increasing sequence from $LAff(K)$. 

For positive elements $a,b\in A$ we write that $a\approx b$ if there is $x\in A$ such that $x^*x=a$ and $xx^*=b$. The relation $\approx$ is an equivalence relation, and it is known that $a\approx b$ implies $a\sim b$.\\

\section{Some Useful Results}
The first part of our paper is based on the work of \cite{1} and will be using some of its lemmas. We will first restate these lemmas as follows for completeness and future references.

For each $\eta>0$ we define a continuous map $f_\eta:\mathbb{R}_+\to[0,1]$ by the following formula:
\[
f_\eta(t) =
\begin{cases}
t/\eta, & 0<t<\eta\\
1, & t\geq \eta.
\end{cases}
\]

\begin{lem}(See [8, Lemma 3.1].) \label{1a} 
Let $A$ be a unital C*-algebra, with $T(A)\neq \emptyset$ and let $a \in M_{k}(A)$ be positive. Suppose that there are $0<\alpha<\beta$ such that $\alpha<d_\tau(a)<\beta$ for every $\tau$ in a closed subset $Y$ of $T(A)$. Then there exists $\varepsilon>0$ and an open nieghborhood $U$ of $Y$, with the property that
\[
\alpha < d_{\tau}((a-\varepsilon)_{+}) < \beta, \forall\tau\in U.
\]
\end{lem}

\begin{lem} (See [8, Lemma 3.2].)\label{1b} 
Let $A$ be a unital separable C*-algebra with nonempty tracial state space, and let $Y\subset T(A)$ be closed. Suppose that $a \in M_{k}(A)$ is a positive element with the property that
\[
\beta-\alpha < d_{\tau}(a) \leq \beta, \forall\tau\in Y
\]
for some $0<\alpha<\beta$. Then there exists $\eta>0$ such that
\[
k-\beta \leq d_\tau(1_k-f_\eta(a)) <k-\beta+2\alpha, \forall\tau\in Y
\]

\end{lem}

\begin{lem}(See [8, Lemma 4.1].)\label{1c}
Let $A$ be a unital simple separable infinite-dimensional C*-algebra and $\tau$ a normalized trace on $A$. Let $0<s<r$ be given. It follows that there are an open neighborhood $U$ of $\tau$ in $T(A)$ and a positive element a in some $M_k(A)$ such that
\[
s<d_\gamma(a)<r, \forall \gamma\in U\
\]
\end{lem}

\begin{lem}(See [8, Lemma 4.2].)\label{1d}
Let $A$ be a unital C*-algebra and $\tau$ a normalized trace on $A$. Let $x$,$y$ be positive elements in $M_k(A)$. Then $d_\tau(y^{1/2}xy^{1/2})\geq d_\tau(x)-d_\tau(1_k-y)$, where $1_k$ denotes the unit of $M_k(A)$.
\end{lem}

The next theorem is from Lin's paper \cite{2} based on work of Cuntz and Pedersen:
\begin{thm}(See [14, Theorem 9.3].)\label{1e}
Let $A$ be a unital simple C*-algebra with nonempty tracial state space, and let $f$ be a strictly positive affine continuous function on $T(A)$. It follows that for any $\varepsilon>0$ there is a positive element $a$ of $A$ such that $f(\tau)=\tau(a), \forall \tau\in T(A)$, and $\left\| a\right\|<\left\| f\right\|+\varepsilon$.
\end{thm}

\section{Rank Functions}

If $a$ is a positive element in $A$ and $\tau\in T(A)$, we denote $\nu_\tau$ the measure induced on the spectrum $\sigma(a)$ of $a$ by $\tau$. Then \mbox{$d_\tau(a)=\nu_\tau((0,\infty)\cap\sigma(a))$} and more generally
\[
d_\tau(f(a))=\nu_\tau(\{t\in\sigma(a):f(t)>0\})
\]
for all nonnegative continuous functions $f$ defined on $\sigma(a)$. (See \cite{1}, under Definition 2.1)

\begin{lem} \label{2a} 
Let $A$ be a unital C*-algebra, with $T(A)\neq \emptyset$ and let $a \in M_{k}(A)$ be positive. Suppose that $Y$ is a compact subset of $T(A)$. Then $\forall  \delta >0$, $\forall \delta'>0$, $\exists \varepsilon >0$, and an open nieghborhood $U$ of $Y$, such that $\forall \gamma \in U$, $\exists \tau, \tau' \in Y$, $dist(\tau, \gamma)<\delta'$, $dist(\tau', \gamma)<\delta'$and
\begin{equation*}
d_{\tau}(a)-\delta < d_{\gamma}((a-\varepsilon)_{+}) < d_{\tau'}(a)+\delta
\end{equation*}
\end{lem}

\begin{proof}
Since $d_{\tau}((a-\varepsilon)_{+}) \nearrow d_{\tau}(a)$ as $\varepsilon \searrow 0$ for each $\tau$, we can fix $\varepsilon_{\tau} >0$, such that $d_{\tau}((a-\varepsilon_{\tau})_{+}) > d_{\tau}(a)-\delta$. \\
Since $\gamma \mapsto d_{\gamma}((a-\varepsilon_{\tau})_{+})$ is lower semi-continuous, we can find an open neighborhood $V_{\tau}$ of $\tau$, such that
\begin{align*}
&dist(\tau, \gamma)<\delta', \forall \gamma \in V_{\tau}\\
&d_{\gamma}((a-\varepsilon_{\tau})_{+})> d_{\tau}(a)-\delta, \forall \gamma \in V_{\tau}
\end{align*}
The family $\{V_{\tau}\}_{\tau \in Y}$ is an open cover of $Y$. Since $Y$ is compact, $Y \subset V_{\tau_{1}} \cup ...\cup V_{\tau_{n}}$ for some $\tau_{1}, ...,\tau_{n} \in Y$.
Set $\varepsilon := \min\{\varepsilon_{\tau_{1}}, ...,\varepsilon_{\tau_{n}}\}$ and $V:=V_{\tau_{1}} \cup ...\cup V_{\tau_{n}}$, so that for each $\gamma \in V$, we have $\gamma \in V_{\tau_{i}}$ for some $i$, and \\
\begin{equation*}
d_{\gamma}((a-\varepsilon)_{+})\geq d_{\gamma}((a-\varepsilon_{\tau_{i}})_{+}) > d_{\tau_i}(a)-\delta
\end{equation*}
On the other hand, for any $\tau'\in T(A)$, let $\nu_{\tau'}$ be the measure induced on $\sigma(a)$ by $\tau'$, we also have
\begin{align*}
d_{\tau'}((a-\varepsilon)_{+})&=\nu_{\tau'}((\varepsilon, +\infty)\cap\sigma(a))\\
&\leq \nu_{\tau'}([\varepsilon, +\infty)\cap\sigma(a))\\
&\leq \nu_{\tau'}((0, +\infty)\cap\sigma(a))\\
&\leq d_{\tau'}(a).
\end{align*}
In particular, 
\[
d_{\tau'}((a-\varepsilon)_{+})\leq\nu_{\tau'}([\varepsilon, +\infty)\cap\sigma(a))<d_{\tau'}(a)+\delta
\]
for all $\tau' \in Y$. By the Portmanteau theorem(\cite{27}, Theorem 13.16), the map $\gamma \mapsto \nu_{\gamma}([\varepsilon, +\infty)\cap\sigma(a))$ is upper semi-continuous, and so the set $W_{\tau'}=\{\gamma\in T(A):\nu_{\gamma}([\varepsilon, +\infty)\cap\sigma(a))<d_{\tau'}(a)+\delta\}$ is open and contains $\tau'$. Without loss of generality, we can assume that the diameter of $W_\tau'$ is less than $\delta'$. Hence, for any $\gamma \in W_{\tau'}$, we have $dist(\tau', \gamma)<\delta'$ and $d_{\gamma}((a-\varepsilon)_{+})<d_{\tau'}(a)+\delta.$\\
Set $W=\bigcup_{\tau'\in Y}W_{\tau'}$ and let $U=W\cap V$. $U$ is an open neighborhood of Y. For any $\gamma\in U, \gamma\in V_{\tau}$ for some $\tau$, and $\gamma\in W_{\tau'}$ for some $\tau'$, hence $dist(\tau, \gamma)<\delta'$, $dist(\tau', \gamma)<\delta'$ and
\begin{equation*}
d_{\tau}(a)-\delta<d_{\gamma}((a-\varepsilon)_{+})<d_{\tau'}(a)+\delta.
\end{equation*}
holds.
\end{proof}

The next lemma is a generalization of Lemma 4.4 of \cite{1}, which generates the \enquote{indicator rank functions}.
\begin{lem} \label{2b}
Let $A$ be a separable unital simple infinite-dimensional C*-algebra whose tracial state space $T(A)$ is nonempty and \mbox{$\partial_e T(A)=X$} is a $F_\sigma$ set. It follows that for any $\delta >0$, and compact subset $Y\subset X$, there is a nonzero positive element $a$ of $A$ with the property that
\begin{align*}
d_{\tau}(a)&<\delta, \forall \tau\in Y,\\
d_{\tau}(a)&=1, \forall\tau\in X\backslash Y.
\end{align*}
\end{lem}

\begin{proof}
Let $\delta$ and $Y$ be given. If $X$ is compact, then this result has already been established in Lemma 4.4 of \cite{1}. Assume now that $X$ is non-compact, so $Y\neq X$. Fix a decreasing sequence $\{U_{n}\}_{n=2}^{\infty}$ of open subsets of $X$ with the property that $Y=\bigcap_{n=2}^{\infty}U_{n}$ and $U_{n}^{c}\neq\emptyset$.\\
Since $X$ is a $F_{\sigma}$ set, there exists an increasing sequence $\{F_{n}\}_{n=2}^{\infty}$ of compact subsets of $X$ such that $X=\bigcup_{n=2}^{\infty}F_{n}$ and $F_{n}=Y\cup E_{n}$, where $\{E_{n}\}_{n=2}^{\infty}$ is an increasing sequence of compact subsets and $E_{n}\subset U_{n}^{c}$ for each $n$. Then by Corollary 11.15 in \cite{3}, since $E_{n}$ and $Y$ are compact, we can use Theorem \ref{1e}  to produce a sequence of $(b_{n})_{n=2}^{\infty}$ in $A_{+}$ with the following properties:
\begin{align*}
\tau(b_{n})&>1-1/n, \forall \tau\in E_{n} \\
\tau(b_{n})&<\delta/(2^{n}n), \forall \tau\in Y \\
\|b_{n}\|&\leq 1.
\end{align*}
For any $\tau\in E_n$ we have
\[
d_{\tau}((b_{n}-1/n)_{+})\geq\tau((b_{n}-1/n)_{+})\geq\tau(b_{n})-1/n>1-2/n.
\]
In particular $(b_{n}-1/n)_{+}\neq 0$. Moreover, for any $\tau\in Y$ we have
\[
d_{\tau}((b_{n}-1/n)_{+})=n\int{(1/n)\chi_{(1/n,\infty)}d\mu_{\tau}}\leq n\tau(b_{n})<\delta/2^{n}.
\]
Set $c_{n}:=2^{-n}(b_{n}-1/n)_{+}$, so that $d_{\tau}(c_{n})>1-2/n$ for each $\tau\in E_n$, $d_{\tau}(c_{n})<\delta/2^{n}$ for each $\tau\in Y$ and $\|c_{n}\|\leq2^{-n}$. Now set
\[
a:=\sum_{n=2}^\infty c_{n}\in A_{+}
\]
If $\tau\in Y$, then using the lower semi-continuity of $d_{\tau}$, we have
\begin{equation*}
\begin{split}
d_{\tau}(a)& \leq\liminf_{k} d_{\tau}(\sum_{n=2}^k c_{n})\\
&\leq\liminf_{k}\sum_{n=2}^k d_\tau(c_{n})\\
&\leq\delta
\end{split}
\end{equation*}
If $\tau\in X\backslash Y$, then $\tau\in E_k$ for all $k$ sufficiently large. It follows that for these same $k$,
\[
d_\tau(a)=d_\tau(\sum_{n=2}^\infty c_n)\geq d_\tau(c_k)\geq 1-2/k
\]
We conclude that $d_\tau(a)\geq 1$ for each such $\tau$. On the other hand, $a\in A$, so $d_\tau(a)\leq1$ for any $\tau\in T(A)$.
\end{proof}

\begin{lem}\label{2c}
Let $A$ be a separable unital simple infinite-dimensional C*-algebra. Suppose that $\partial_eT(A)=X$ is a nonempty $F_\sigma$ set. If there exist a compact subset $F\subset X$ and some $0<\varepsilon<1$ such that $\mu_\tau(F)>\varepsilon, \forall\tau\in\partial X$, then $\partial X\cup F$ is a compact subset of $T(A)$.
\end{lem}

\begin{proof}
Since there exist a compact subset $F\subset X$ and some $0<\varepsilon<1$ such that $\mu_\tau(F)>\varepsilon, \forall\tau\in\partial X$, then by previous lemma, there exists a nonzero positive element $a$ of $A$ with the property that
\begin{align*}
d_{\tau}(a)&<\varepsilon/2, \forall \tau\in F,\\
d_{\tau}(a)&=1, \forall\tau\in X\backslash F.
\end{align*}
Hence for each $\tau\in\partial X$,
\begin{align*}
d_{\tau}(a)&=\int_X d_{\gamma}(a)d\mu_\tau\\
&=\int_F d_{\gamma}(a)d\mu_\tau+\int_{X\backslash F} d_{\gamma}(a)d\mu_\tau\\
&\leq \varepsilon/2\cdot\mu_\tau(F)+\mu_\tau(X\backslash F)\\
&\leq 1-\varepsilon/2
\end{align*}
Since $d_{\tau}(a)$ is lower semi-continuous, $T=\{\tau\in T(A)|d_{\tau}(a)\leq 1-\varepsilon/2\}$ is compact. By the compactness of $\overline{X}$, $T\cap \overline{X}=\partial X\cup F$ is compact.
\end{proof}

\begin{rem}
This lemma plays an important role in the proof of Theorem 5.1. We can verify that for a general extreme boundary $X$, if $X$ has the property that there exists some compact subset $F\subset X$ such that $\partial X\cup F$ is compact, then $X$ must be a $F_\sigma$ set. To see this, consider $\overline{X}\backslash (\partial X\cup F)$ which is a $F_\sigma$ set. Then $X=\overline{X}\backslash\partial X$ is also a $F_\sigma$ set.

\end{rem}

The next lemma is a generalization of Lemma 4.5 of \cite{1}.
\begin{lem} \label{2d}
Let $A$ be a unital simple separable infinite dimensional C*-algebra. Suppose that $X=\partial_e T(A)$ is a nonempty $F_\sigma$ set. Let $a\in M_N(A)$ be positive, and let there be a given compact subset Y of X and $\delta>0$. It follows that there is a positive element $b$ of $M_N(A)$ with the following properties:
\begin{align*}
d_\tau(b)&=d_\tau(a), \forall \tau\in X\backslash Y\\
d_\tau(b)&\leq\delta, \forall \tau\in Y
\end{align*}
\end{lem}

\begin{proof}
By possibly replacing $A$ with $M_N(A)$, it suffices to prove this lemma in the case of $N=1$. Use Lemma \ref{2b} to find a positive element $h$ of $A$ satisfying
\begin{align*}
d_\tau(h)&<\delta, \forall \tau\in Y\\
d_\tau(h)&=1, \forall \tau\in X\backslash Y.
\end{align*}
Since $X$ is a nonempty $F_\sigma$ set, there exists a sequence of open subsets of $X$, $V_1\subset V_2\subset V_3\subset...$, such that $\overline{V_i}\subset X\backslash Y$ for each $i$ and \mbox{$\bigcup_{i=1}^\infty V_i=X\backslash Y$}. Now the rest of the proof is similar to that of \mbox{Lemma 4.5} of \cite{1}. Trivially,
\[
1-1/2i<d_\tau(h)\leq 1, \forall\tau\in\overline{V_i},
\]
and so Lemma \ref{1b} applied for $k=\beta=1$ and $\alpha=1/2i$ yields $\eta_i>0$ such that
\[
d_\tau(1-f_{\eta_i}(h))<1/i, \forall\tau\in\overline{V_i}
\]
where $f_{\eta}$ is defined in Section 3. To simplify notation in the remainder of the proof, relabel $f_{\eta_i}(h)$ as $h_i$. We may assume that the sequence $(\eta_i)$ is decreasing so that the sequence $(h_i)$ is increasing. Since $d_\tau(h)=d_\tau(f_\eta(h))$ for any $\tau\in T(A)$ and $\eta>0$, it follows that
\begin{align*}
d_\tau(h_i)&\leq\delta, \forall \tau\in Y\\
d_\tau(h_i)&=1, \forall \tau\in X\backslash Y.
\end{align*}
Set $a_i:=a^{1/2}h_ia^{1/2}$. Since $a_i=a^{1/2}h_ia^{1/2}\approx h_i^{1/2}ah_i^{1/2}\precsim a$, we have
\[
d_\tau(a_i)\leq d_\tau(a), \forall\tau\in X\backslash Y.
\]
Also, since $a_i=a^{1/2}h_ia^{1/2}\precsim h_i$, we have
\[
d_\tau(a_i)\leq d_\tau(h_i)<\delta, \forall\tau\in Y.
\]
For our lower bound, we observe that by Lemma \ref{1d} we have for any $\tau\in\overline{V_i}$:
\[
d_\tau(a_i)=d_\tau(h_i^{1/2}ah_i^{1/2})\geq d_\tau(a)-d_\tau(1-h_i)>d_\tau(a)-1/i.
\]
Therefore we have
\begin{align*}
&d_\tau(a_i)<\delta, \forall \tau\in Y\\
&d_\tau(a)-1/i<d_\tau(a_i)<d_\tau(a), \forall \tau\in \overline{V_i}.
\end{align*}
Since $h_i\leq h_{i+1}$ and $a_i\approx h_i^{1/2}ah_i^{1/2}\precsim h_{i+1}^{1/2}ah_{i+1}^{1/2}\approx a_{i+1}$, the increasing sequence $(\langle a_i\rangle)_{i=1}^\infty$ has a supremum $y$, where $y=\langle b\rangle$ for some positive element $b$ of $A\otimes\mathcal{K}$ by \cite{53}. Since each $d_\tau$ is a supremum preserving state on $Cu(A)$, we conclude that
\begin{align*}
d_\tau(b)&<\delta, \forall \tau\in Y\\
d_\tau(b)&=d_\tau(a), \forall \tau\in X\backslash Y.
\end{align*}
as desired.
\end{proof}

The next lemma is a generalization of Lemma 4.6 of \cite{1}. Note that Lemma 4.6 of \cite{1} uses Lemma 4.4 and 4.5 of \cite{1}, which we have generalized in Lemma \ref{2b} and \ref{2d} in this paper. Therefore, by adapting the proof of Lemma 4.6 of \cite{1}, we know the following lemma is true. 
\begin{lem}\label{2h}
Let $A$ be a unital simple separable C*-algebra with strict comparison of positive elements and $T(A)$ is a nonempty $F_\sigma$ set. $F$ is a compact subset of the extreme boundary $X=\partial_e T(A)$ and $a,b\in A$ is positive. Suppose that there are $0<\alpha<\beta<\gamma\leq 1$ and sets $U,V\subset F$ that are relatively open in $F$ with the property that $\alpha<d_\tau(a)<\beta, \forall\tau\in U$ and $\beta<d_\tau(b)<\gamma, \forall\tau\in V$. It follows that for any closed set $K\subset U\cup V$, there is a positive element $c$ of $M_2(A)$ with the property that
\[
\alpha<d_\tau(c)<\gamma, \forall\tau\in K.
\]
\end{lem}

By using Lemma \ref{2d} and adapting the proof of Theorem 5.2 of \cite{1}, we obtain a generalization of it.
\begin{lem} \label{2g}
Let $A$ be a unital simple separable C*-algebra. Suppose that $X=\partial_e T(A)$ is a nonempty $F_\sigma$ set. Assume further that for each $m\in \mathbb{N}$ and any compact subset $F\subset \partial_e T(A)=X$, there is $x\in Cu(A)$ with the property that $md_\tau(x)\leq1\leq(m+1)d_\tau(x),\forall\tau\in$ F.
It follows that for any $f\in \mathit{Aff}(T(A))$, $\varepsilon>0$, any compact subset $F\subset X$, there is positive $h\in A\otimes \mathcal{K}$, such that
\[
|d_\tau(h)-f(\tau)|<\varepsilon,\forall\tau\in F.
\]
\end{lem}

Similarly, by using Lemma \ref{2g}, \ref{1a}, \ref{2d} and \ref{2h}, and by adapting the proof of Theorem 5.4 of \cite{1} the next lemma can be extracted from an intermediate result in the proof of Theorem 5.4 of \cite{1}.

\begin{lem}\label{2e}
Let $A$ be a unital simple separable C*-algebra. Suppose that the extreme boundary $X$ of $T(A)$ is nonempty, $F_\sigma$ and of finite covering dimension. Then for any compact subset $F\subset X$, and for each $0\leq r'<r<1$, there exists a positive element a in some $M_N(A)$ with the property that
\[
r'<d_\tau(a)<r, \forall\tau\in F.
\]
\end{lem}

\begin{lem} \label{2i}
Let $A$ be a unital simple separable C*-algebra with strict comparison of positive elements. Suppose that the extreme boundary $X$ of $T(A)$ is nonempty, zero dimensional, $F_\sigma$ and has the tightness property. It follows that for any $f\in \mathit{Aff}(T(A))$ with $\left\| f\right\|=1$, $0<\varepsilon<1$, any compact subset $F\subset X$, there is $h\in (A\otimes \mathcal{K})_+$, such that
\begin{align*}
&|d_\tau(h)-f(\tau)|<\varepsilon,\forall\tau\in F\\
&d_\tau(h)\leq 2, \forall\tau\in X\backslash F.
\end{align*}
\end{lem}

\begin{proof}
First of all, since $A$ is simple and has strict comparison of positive elements, by Theorem 4.4.1 in \cite{13}, for any element $\langle a\rangle\in Cu(A)$ of $A$, if its rank function is bounded, that is, if there exists some $m>0$ such that \mbox{$\forall\tau\in T(A)$}, $d_\tau(a)\leq m$, then we can always find an element $a'$ in some $M_N(A)$ such that $a$ and $a'$ are Murray-von Neumann equivalent. Because we only concern about rank functions, and the elements we are dealing with all satisfy the condition that their rank functions are bounded, in this proof we will no longer distinguish between elements in $A\otimes \mathcal{K}$ and elements in the matrix algebras of $A$.\\
Since $X$ is zero dimensional, we can find a compact subset $F_1$ of $X$ containing $F$ such that $X\backslash F_1$ is closed. By Lemma \ref{2e} and Lemma \ref{2g}, there exists $h_1\in (A\otimes \mathcal{K})_+$ satisfying
\[
|d_\tau(h_1)-f(\tau)|<\varepsilon/2,\forall\tau\in F_1.
\]
There exists some $m>0$ such that $d_\tau(h_1)\leq m, \forall\tau\in T(A)$. Since $X$ has the tightness property, we can find a compact subset $F_2$ of $X$ such that $\mu_\tau(F_2)>1-\varepsilon_1, \forall\tau\in\partial X$ for $\varepsilon_1=(1-\varepsilon)/2m$. Since $X\backslash F_1$ is closed, $F_2\backslash F_1$ is compact. Then by Lemma \ref{2d} there exists $h_2\in (A\otimes \mathcal{K})_+$ such that
\begin{align*}
&|d_\tau(h_2)-f(\tau)|<\varepsilon/2,\forall\tau\in F_1\\
&d_\tau(h_2)\leq \delta_1, \forall\tau\in F_2\backslash F_1
\end{align*}
for $\delta_1=(1-\varepsilon)/2$.\\
Since $d_\tau(h_2)$ is affine lower semi-continuous, $d_\tau(h_2)=\int_{X} d_\gamma(h_2) d\mu_\tau$ for each $\tau\in\partial X$. Then,
\begin{align*}
|d_\tau(h_2)| &\leq | \int_{F_1} d_\gamma(h_2) d\mu_\tau|+| \int_{F_2\backslash F_1} d_\gamma(h_2) d\mu_\tau|+| \int_{X\backslash F_2} d_\gamma(h_2) d\mu_\tau| \\
&\leq(1+\varepsilon/2) \cdot\mu_\tau(F_1)+\delta_1\cdot\mu_\tau(F_2\backslash F_1)+m\cdot\mu_\tau(X\backslash F_2)\\
&\leq1+\varepsilon/2+\delta_1+m\cdot\varepsilon_1\\
&< 2
\end{align*}
for each $\tau\in\partial X$.\\
Since $f$ is uniformly continuous on $T(A)$, there exists $\delta_2>0$ such that for any $\tau, \gamma \in T(A)$, $|f(\tau)-f(\gamma)|<\varepsilon/6$ holds whenever \mbox{$dist(\tau, \gamma)<\delta_2$}. Since $F_1$ is compact, by Lemma \ref{2a}, we obtain $\varepsilon_2>0$ and an open neighborhood $V$ of $F_1$, such that for any $\gamma\in V, \exists\tau, \tau'\in F_1$, such that \mbox{$dist(\tau, \gamma)<\delta_2$}, $dist(\tau', \gamma)<\delta_2$ and
\[
d_\tau(h_2)-\varepsilon/6<d_\gamma((h_2-\varepsilon_2)_+)<d_{\tau'}(h_2)+\varepsilon/6.
\]
Denote $h_3:=(h_2-\varepsilon_2)_+$. Then
\begin{align*}
d_\gamma(h_3)-f(\gamma) &= d_\gamma(h_3)-d_\tau(h_2)+d_\tau(h_2)-f(\tau)+f(\tau)-f(\gamma)\\
&\geq -\varepsilon/6-\varepsilon/2-\varepsilon/6\\
&\geq-5\varepsilon/6
\end{align*}
and
\begin{align*}
d_\gamma(h_3)-f(\gamma) &= d_\gamma(h_3)-d_{\tau'}(h_2)+d_{\tau'}(h_2)-f(\tau')+f(\tau')-f(\gamma)\\
&\leq \varepsilon/6+\varepsilon/2+\varepsilon/6\\
&\leq 5\varepsilon/6
\end{align*}
for all $\gamma\in V$.\\
Let $\nu_\tau$ be the measure induced on $\sigma(h_2)$ by $\tau$. Obviously we have $d_\tau(h_3)=d_\tau((h_2-\varepsilon_2)_+)\leq\nu_\tau([\varepsilon_2,\infty)\cap\sigma(h_2))\leq d_\tau(h_2)<2$ for all $\tau\in\partial X\cup F_1$. By the Portmanteau Theorem(\cite{27}, Theorem 13.16), the map $\gamma \mapsto \nu_{\gamma}([\varepsilon_2, +\infty)\cap\sigma(h_2))$ is upper semi-continuous, and so the set $W=\{\gamma\in T(A):\nu_{\gamma}([\varepsilon_2, +\infty)\cap\sigma(h_2))<2\}$ is open and contains $\partial X\cup F_1$. Moreover, for any $\gamma \in W$, we have $d_{\gamma}(h_3)<2$.\\
Since $X\backslash W$ is a compact subset of $X$, by Lemma \ref{2d} there exists $h_4\in (A\otimes \mathcal{K})_+$ such that
\begin{align*}
&|d_\tau(h_4)-f(\tau)|<\varepsilon,\forall\tau\in W\\
&d_\tau(h_4)\leq 2, \forall\tau\in X\backslash W.
\end{align*}
This proves the lemma since $F\subset W$, $\varepsilon<1$ and $\left\| f\right\|=1$.
\end{proof}

\section{Rank Functions on Zero-dimensional Extreme Boundaries}
\begin{thm}\label{3a}
Let $A$ be a unital simple separable C*-algebra with strict comparison of positive elements. Suppose further that the extreme boundary $X$ of $T(A)$ is nonempty, zero dimensional, $F_\sigma$ and has the tightness property.
It follows that $\forall f\in \mathit{SAff}(T(A))$, $\exists h\in (A\otimes \mathcal{K})_+$ such that $d_\tau(h)=f(\tau)$ for each $\tau\in T(A)$.\\
In particular, this theorem holds if $X$ is countable and has finitely many cluster points. 
\end{thm}

\begin{proof}
We will first prove that $\forall f\in \mathit{Aff}(T(A))$, $\forall$ $0<\varepsilon<1$, \mbox{$\exists h\in (A\otimes \mathcal{K})_+$}, such that
\begin{equation*}
|d_\tau(h)-f(\tau)|<\varepsilon, \forall \tau\in T(A). \tag{$\ast$}
\end{equation*}
We only need to establish ($\ast$) on $\partial_e T(A)=X$. \\
For convenience, we say a subset $Y\subset X$ is clopen if both $Y$ and $\overline{X}\backslash Y$ are compact. In addition, for the same reason given in the proof of Lemma 4.8, in this proof we will no longer distinguish between elements in $A\otimes \mathcal{K}$ and elements in the matrix algebras of $A$.\\
Let $f\in \mathit{Aff}(T(A))$ and $\varepsilon>0$ be given. Without loss of generality assume $\left\| f\right\|=1$. Since $f$ is uniformly continuous on $T(A)$, there exists $\delta'>0$ such that for any $\tau, \gamma \in T(A)$, $|f(\tau)-f(\gamma)|<\varepsilon/6$ holds whenever $dist(\tau, \gamma)<\delta'$. \\
Since $X$ has the tightness property, there is a compact subset $F\subset X$ such that  $\mu_\tau(F)>1-\varepsilon/48, \forall\tau\in\partial X$. By Lemma \ref{2i}, there exists $h_1\in (A\otimes \mathcal{K})_+$, such that
\begin{align*}
&|d_\tau(h_1)-f(\tau)|<\varepsilon/6, \forall \tau\in F\\
&|d_\tau(h_1)|\leq 2, \forall \tau\in X
\end{align*}
Since $d_\tau(h_1)$ is affine lower-semicontinuous, $d_\tau(h_1)=\int_{X} d_\gamma(h_1) d\mu_\tau$ for each $\tau\in\partial X$. Similar equation holds for $f$ since it is affine continuous. Then,
\begin{align*}
|d_\tau(h_1)-f(\tau)| &\leq | \int_{F} d_\gamma(h_1) d\mu_\tau-\int_{F}f d\mu_\tau|+|\int_{X\backslash F} d_\gamma(h_1) d\mu_\tau-\int_{X\backslash F}f d\mu_\tau| \\
&\leq\varepsilon/6 \cdot\mu_\tau(F)+3\cdot\mu_\tau(X\backslash F)\\
&\leq\varepsilon/6+3\varepsilon/48\\
&<\varepsilon/2
\end{align*}
for each $\tau\in\partial X$.\\
By Lemma \ref{2c}, $\partial X\cup F$ is compact. Then by Lemma \ref{2a}, we obtain $\varepsilon_1>0$ and an open neighborhood $V$ of $\partial X\cup F$. For any $\gamma\in V$, \mbox{$\exists\tau, \tau'\in\partial X\cup F$}, such that $dist(\tau, \gamma)<\delta'$, $dist(\tau', \gamma)<\delta'$ and
\[
d_\tau(h_1)-\varepsilon/6<d_\gamma((h_1-\varepsilon_1)_+)<d_{\tau'}(h_1)+\varepsilon/6.
\]
Denote $h_2:=(h_1-\varepsilon_1)_+$. Then
\begin{align*}
d_\gamma(h_2)-f(\gamma) &= d_\gamma(h_2)-d_\tau(h_1)+d_\tau(h_1)-f(\tau)+f(\tau)-f(\gamma)\\
&\geq -\varepsilon/6-\varepsilon/2-\varepsilon/6\\
&\geq-5\varepsilon/6
\end{align*}
and
\begin{align*}
d_\gamma(h_2)-f(\gamma) &= d_\gamma(h_2)-d_{\tau'}(h_1)+d_{\tau'}(h_1)-f(\tau')+f(\tau')-f(\gamma)\\
&\leq -\varepsilon/6+\varepsilon/2+\varepsilon/6\\
&\leq 5\varepsilon/6
\end{align*}
for all $\gamma\in V$.\\
Since $V$ is an open neighborhood of the compact subset $\partial X\cup F$, and since $X$ is zero-dimensional, we can replace $V$ by a open set $V'$ such that $\partial X\cup F\subset V'\subset V$ and $\overline{X}\backslash V'$ is clopen. Denote $Y=\overline{X}\backslash V'$. Then by Lemma \ref{2d}, there exists $h_3\in (A\otimes \mathcal{K})_+$, such that
\begin{align*}
d_\tau(h_3)&=d_\tau(h_2), \forall\tau\in X\backslash Y\\
d_\tau(h_3)&\leq\varepsilon/6, \forall\tau\in Y.
\end{align*}
So
\begin{align*}
&|d_\tau(h_3)-f(\tau)|<5\varepsilon/6, \forall\tau\in X\backslash Y\\
&d_\tau(h_3)\leq\varepsilon/6, \forall\tau\in Y.
\end{align*}
Then it suffices to find $l\in (A\otimes \mathcal{K})_+$ satisfying
\begin{align*}
&d_\tau(l)\leq\varepsilon/6, \forall \tau \in X\backslash Y\\
&|d_\tau(l)-f(\tau)|<\varepsilon/3, \forall \tau \in Y.
\end{align*}
By setting $h:=h_3\oplus l$, we get
\[
|d_\tau(h)-f(\tau)|\leq\varepsilon, \forall \tau\in X
\]
as desired.
The construction of such $l$ is as follows. 

Since $Y$ is clopen, by Lemma \ref{2b}, there exists $t\in A$ with the property $d_\tau(t)<\varepsilon/24$ on $Y$ and $d_\tau(t)=1$ on $X\backslash Y$. For $\tau\in\partial X$, \mbox{$d_\tau(t)\geq1\cdot\mu_\tau(X\backslash Y)\geq 1-\varepsilon/48$}. So we have
\begin{align*}
&d_\tau(t)<\varepsilon/24, \forall \tau \in Y\\
&1-\varepsilon/48\leq d_\tau(t)\leq 1, \forall \tau \in\overline{X}\backslash Y.
\end{align*}
Since $\overline{X}\backslash Y$ is compact, by Lemma \ref{1a}, there are $\eta>0$, such that
\begin{align*}
d_\tau((t-\eta)_+)>1-\varepsilon/48, \forall\tau\in \overline{X}\backslash Y.
\end{align*}
Denote $\nu_\tau$ the measure induced on the spectrum $\sigma(t)$ of $t$ by $\tau$. Since $d_\tau((t-\eta)_+)=\nu_\tau((\eta,\infty)\cap\sigma(t))$, then
\begin{align*}
\nu_\tau((0,\eta]\cap\sigma(t))&=d_\tau(t)-\nu_\tau((\eta,\infty)\cap\sigma(t))\\
&<1-(1-\varepsilon/48)\\
&\leq\varepsilon/48, \forall\tau\in \overline{X}\backslash Y.
\end{align*}
Then $\forall\tau\in \overline{X}\backslash Y$
\begin{align*}
d_\tau(1-f_\eta(t))&=\nu_\tau([0,\eta)\cap\sigma(t))\\
&=\nu_\tau((0,\eta)\cap\sigma(t))+\nu_\tau(\{0\}\cap\sigma(t))\\
&=\nu_\tau((0,\eta)\cap\sigma(t))+1-d_\tau(t)\\
&<\varepsilon/48+1-(1-\varepsilon/48)\\
&\leq\varepsilon/24
\end{align*}
where $f_\eta$ is defined in Section 3. On the other hand, $\forall\tau\in Y$
\begin{align*}
d_\tau(1-f_\eta(t))&=\nu_\tau((0,\eta)\cap\sigma(t))+1-d_\tau(t)\\
&\geq1-d_\tau(t)\\
&>1-\varepsilon/24.
\end{align*}
Now, since $Y$ is compact, by Lemma \ref{2i}, we can find $x \in A\otimes \mathcal{K}$ such that
\begin{align*}
&|d_\tau(x)-f(\tau)|<\varepsilon/6, \forall \tau \in Y\\
&d_\tau(x)\leq 2, \forall\tau\in X\backslash Y.
\end{align*}
Hence, by Theorem 4.4.1 in \cite{13} we can find $x'\in M_4(A)$ which is Murray-von Neumann equivalent to $x$.
Let \mbox{$t'=\oplus_{j=1}^4 (1-f_\eta(t))\in M_4(A)$}. Consider $l=t'^{1/2}x't'^{1/2}$.\\
Since $l\precsim x'$, we have $d_\tau(l)\leq d_\tau(x')<f(\tau)+\varepsilon/6$ on $Y$. On the other hand, by Lemma 3.4,
\begin{align*}
d_\tau(l)&\geq d_\tau(x')-d_\tau(1_4-t')\\
&\geq d_\tau(x')-4\cdot d_\tau(f_\eta(t))\\
&\geq d_\tau(x')-4\cdot d_\tau(t)\\
&>f(\tau)-\varepsilon/6-4\cdot \varepsilon/24\\
&>f(\tau)-\varepsilon/3
\end{align*}
for any $\tau\in Y$. For $\tau\in X\backslash Y$, $d_\tau(l)\leq d_\tau(t')= 4\cdot d_\tau(1-f_\eta(t))<\varepsilon/6$.\\
This finishes the proof of ($\ast$).\\
Now suppose that $A$ has strict comparison. The final conclusion of the Theorem then follows from the proof of Theorem 2.5 of \cite{8}, which shows how one produces an arbitrary $f\in \mathit{SAff}(T(A))$ by taking suprema.
\end{proof}

\section{Proof of the Main Result}
In this section, we will prove Theorem \ref{4a}. First recall that $T_\infty (A)$ is the set of all traces on $A_\infty$ induced by the trace $(x_n)^\infty_{n=1}\mapsto \lim_{n\rightarrow\omega}{\tau_n(x_n)}$ on $\ell^\infty(A)$ where $(\tau_n)^\infty_{n=1}$ is a sequence in $T(A)$ and $\omega\in\beta\mathbb{N}\backslash\mathbb{N}$ is a free ultrafilter(See \cite{14}). If we choose the sequence $(\tau_n)^\infty_{n=1}$ in $Y\subset T(A)$ instead, then write $T^Y_\infty (A)$ for the collection of those traces arising in the same fashion. $T^Y_\infty (A)$ is clearly a subset of $T_\infty (A)$.

Based on the definition of the uniformly tracially large cpc order zero map which is given in Definition 2.2 of \cite{14}, we now introduce another definition.
\begin{dfn}
Let $A$ be a separable unital C*-algebra with \mbox{$T(A)\neq\emptyset$}. Let $Y\subset T(A)$ be nonempty. A completely positive and contractive order zero map $\Phi : M_k\rightarrow A_\infty$ is uniformly tracially large on $Y$ if $\tau(\Phi(1_k))=1$ for all $\tau\in T^Y_\infty(A)$.
\end{dfn}

The following lemma can be obtained by adapting the proof of Lemma 2.3 of \cite{14}.

\begin{lem}\label{4y}
Let $A$ be a separable unital C*-algebra with $T(A)\neq\emptyset$ and let $Y\subset T(A)$ be nonempty. Let $\Phi : M_k\rightarrow A_\infty$ is a cpc order zero map. Then $\Phi$ is uniformly tracially large on $Y$ if and only if any lifting $(\phi_n) : M_k\rightarrow \ell^\infty(A)$ of $\Phi$ to a sequence of cpc order zero maps satisfies
\[
\lim_{n\rightarrow\infty}\min_{\tau\in Y}\tau(\phi_n(1_k))=1.
\]
\end{lem}
The idea of the proof of Theorem \ref{4a} is similar to that of Theorem 4.6 of \cite{14}, in which its Lemma 3.5 plays a key role. Note that if the result of Lemma 3.5 of \cite{14} holds when $\partial_eT(A)$ is replaced by a nonempty subset $Y$ of $\partial_eT(A)$, then we can use the same type of argument as in Section 4 of \cite{14} and get similar versions of Lemma 4.1, Lemma 4.2, Lemma 4.3, Proposition 4.4 and Lemma 4.5 of \cite{14} by replacing $T_\infty (A)$ by $T^Y_\infty (A)$. Therefore the following result is true:

\begin{lem}\label{4b}
Let $A$ be a simple separable unital nuclear nonelementary C*-algebra with $T(A)\neq\emptyset$ and let $Y\subset\partial_eT(A)$ be nonempty. Suppose that there exists some $m\in\mathbb{N}$ such that for every finite set $\mathcal{F}\subset A$ and $\varepsilon>0$, there exists cpc order zero maps $\psi^{(0)}, ..., \psi^{(m)}: M_k\rightarrow A$ such that
\[
\left\| [\psi^{(i)}(x), y]\right\|\leq\varepsilon \left\| x\right\|
\]
for all $i\in\{0, ..., m\}$, $x\in M_k$, $y\in\mathcal{F}$ and such that for each $\tau\in Y$, there exists $i(\tau)\in\{0, ..., m\}$ such that $\tau(\psi^{(i(\tau))}(1_k))>1-\varepsilon$. Then there exists a cpc order zero map $\phi^{(0)}: M_k\rightarrow A$ such that
\[
\left\| [\phi^{(0)}(x), y]\right\|\leq\varepsilon \left\| x\right\|
\]
for $x\in M_k$, $y\in\mathcal{F}$ and such that for each $\tau\in Y$, $\tau(\phi^{(0)}(1_k))>1-\varepsilon$.
\end{lem}

Note that Lemma 3.5 of \cite{14} still holds if we remove the assumption of the compactness of the extreme boundary and restrict to a compact subset of it:
\begin{lem}\label{4c}
Let $m\geq 0$, $k\geq 2$ and let $A$ be a simple separable unital nuclear nonelementary C*-algebra with $T(A)\neq\emptyset$. Then for any compact subset $Y\subset\partial_eT(A)$ such that \mbox{dim$(Y)\leq m$}, each finite set $\mathcal{F}\subset A$ and $\varepsilon>0$, there exists cpc order zero maps $\psi^{(0)}, ..., \psi^{(m)}: M_k\rightarrow A$ such that
\[
\left\| [\psi^{(i)}(x), y]\right\|\leq\varepsilon \left\| x\right\|
\]
for all $i\in\{0, ..., m\}$, $x\in M_k$, $y\in\mathcal{F}$ and such that for each $\tau\in Y$, there exists $i(\tau)\in\{0, ..., m\}$ such that $\tau(\psi^{(i(\tau))}(1_k))>1-\varepsilon$.
\end{lem}

The next lemma follows immediately from Lemma \ref{4b} and Lemma \ref{4c}:

\begin{lem}\label{4d}
Let $m\geq 0$, $k\geq 2$ and let $A$ be a simple separable unital nuclear nonelementary C*-algebra with $T(A)\neq\emptyset$. Then for any compact subset $Y\subset\partial_eT(A)$ such that \mbox{dim$(Y)\leq m$}, each finite set $\mathcal{F}\subset A$ and $\varepsilon>0$, there exists a cpc order zero map $\phi^{(0)}: M_k\rightarrow A$ such that
\[
\left\| [\phi^{(0)}(x), y]\right\|\leq\varepsilon \left\| x\right\|
\]
for $x\in M_k$, $y\in\mathcal{F}$ and such that for each $\tau\in Y$, $\tau(\phi^{(0)}(1_k))>1-\varepsilon$.
\end{lem}

Now we can prove Theorem \ref{4a}.

\begin{proof}
The implication of $(1)\Rightarrow(2)$ has already been established by M. R$\o$rdam in Corollary 4.6 of \cite{19} without the assumption on $T(A)$. Now assume $A$ has strict comparison. Since $X$ has the tightness property, there is a compact subset $Y\subset X$, such that for any $\gamma\in\overline {X} \backslash X$, $\mu_\gamma(Y)>1-\varepsilon/2$. Then by Lemma \ref{4d}, there is $\phi^{(0)}:M_k\rightarrow A$ such that
\[
\left\| [\phi^{(0)}(x), y]\right\|\leq\varepsilon \left\| x\right\|
\]
and $\tau(\phi^{(0)}(1_k))>1-\varepsilon/2$ for any $\tau\in Y$.
Define a map $g: T(A)\rightarrow \mathbb{R}$ as follows:
\[
g(\tau)=\tau(\phi^{(0)}(1_k)).
\]
$g$ is affine and continuous. So for $\gamma\in\overline {X} \backslash X$,
\begin{align*}
\gamma(\phi^{(0)}(1_k))&=g(\gamma)\\
&=\int_Xgd\mu_\gamma\geq\int_Ygd\mu_\gamma\\
&\geq(1-\varepsilon/2)(1-\varepsilon/2)\\
&>1-\varepsilon.
\end{align*}
Let $U=\{\tau\in \overline {X}:g(\tau)>1-\varepsilon\}$. $U$ is an open neighborhood of $Y\cup(\overline {X} \backslash X)$. So $\overline {X}\backslash U=X\backslash U$ is compact. Applying Lemma \ref{4c} to $X\backslash U$, we find $\phi^{(1)}:M_k\rightarrow A$ such that
\[
\left\| [\phi^{(1)}(x), y]\right\|\leq\varepsilon \left\| x\right\|
\]
and $\tau(\phi^{(1)}(1_k))>1-\varepsilon$ for any $\tau\in X\backslash U$.
Now we have two cpc order zero maps $\phi^{(0)}, \phi^{(1)}: M_k\rightarrow A$ such that
\[
\left\| [\phi^{(i)}(x), y]\right\|\leq\varepsilon \left\| x\right\|
\]
for $i\in\{0, 1\}$, $x\in M_k$, $y\in\mathcal{F}$ and such that for each $\tau\in X$, there exists $i(\tau)\in\{0, 1\}$ such that $\tau(\phi^{(i(\tau))}(1_k))>1-\varepsilon$.
By Lemma \ref{4b}, we could reduce the number of cpc order zero maps to one. Therefore using the same type of argument as in Theorem 3.6 of \cite{14}, $A$ admits uniformly tracially large cpc order zero maps $M_k\rightarrow A_\infty\cap A'$. Then by Theorem 2.6 of \cite{14}, $A$ is $\mathcal{Z}$-stable.
\end{proof}

We now introduce another definition.
\begin{dfn}
For a unital separable C*-algebra $A$ with $T(A)\neq\emptyset$, denote $X$ the extreme boundary of $T(A)$. From previous discussion we know that there is a Borel probability measure $\mu_\tau$ on $X$ representing each $\tau\in T(A)$.\\ Let $Y, Z\subset X$ be nonempty and denote $\partial Z=\overline {Z} \backslash Z$. We say that $Z$ is weakly tight relative to $Y$ if
\[
\mu_\tau(Y)=1
\]
for all $\tau\in \partial Z$.
\end{dfn}

\begin{lem}\label{4z}
Let $k\geq 2$ and suppose $A$ is a simple nuclear separable unital nonelementary C*-algebra and $T(A)$ is nonempty and of finite covering dimension. Denote $X=\partial_e T(A)$ and let $Y, Z\subset X$ be nonempty. Suppose $Z$ is tight relative to $Y$. If there exists a cpc order zero map $\Phi : M_k\rightarrow A_\infty\cap A'$ which is uniformly tracially large on $Y$, then there exists a cpc order zero map $\Phi' : M_k\rightarrow A_\infty\cap A'$ which is uniformly tracially large on $Y\cup Z$.
\end{lem}

\begin{proof}
We will show that for each finite set $\mathcal{F}\subset A$ and $\varepsilon>0$, there exists a cpc order zero map $\phi: M_k\rightarrow A$ such that
\[
\left\| [\phi(x), y]\right\|\leq\varepsilon \left\| x\right\|
\]
for $x\in M_k$, $y\in\mathcal{F}$ and such that for each $\tau\in Y\cup Z$, \mbox{$\tau(\phi(1_k))>1-\varepsilon$}. Then by taking a nested sequence $(\mathcal{F}_n)^\infty_{n=1}$ of finite subset of $A$ whose union is dense in $A$ and setting $\varepsilon=1/n$ for each $\mathcal{F}_n$, the sequence $(\phi_n)^\infty_{n=1}$ induces a cpc order zero map \mbox{$\Phi' : M_k\rightarrow A_\infty\cap A'$} which is uniformly tracially large on $Y\cup Z$ by Lemma \ref{4y}.

Since there exists a cpc order zero map $\Phi : M_k\rightarrow A_\infty\cap A'$ which is uniformly tracially large on $Y$, by Lemma \ref{4y} there is a cpc order zero map $\phi^{(0)}:M_k\rightarrow A$ such that
\[
\left\| [\phi^{(0)}(x), y]\right\|\leq\varepsilon \left\| x\right\|
\]
and $\tau(\phi^{(0)}(1_k))>1-\varepsilon$ for any $\tau\in Y$.
Define a map $g: T(A)\rightarrow \mathbb{R}$ as follows:
\[
g(\tau)=\tau(\phi^{(0)}(1_k)).
\]
$g$ is affine and continuous. So for $\gamma\in\partial Z$,
\begin{align*}
\gamma(\phi^{(0)}(1_k))&=g(\gamma)\\
&=\int_Xgd\mu_\gamma\geq\int_Ygd\mu_\gamma\\
&>1-\varepsilon.
\end{align*}
Let $U=\{\tau\in \overline {Z}:g(\tau)>1-\varepsilon\}$. $U$ is an open neighborhood of $\partial Z$. So $\overline Z\backslash U=Z\backslash U$ is compact. Applying Lemma \ref{4c} to $Z\backslash U$, we find another cpc order zero map $\phi^{(1)}:M_k\rightarrow A$ such that
\[
\left\| [\phi^{(1)}(x), y]\right\|\leq\varepsilon \left\| x\right\|
\]
and $\tau(\phi^{(1)}(1_k))>1-\varepsilon$ for any $\tau\in Z\backslash U$.
Now we have two cpc order zero maps $\phi^{(0)}, \phi^{(1)}: M_k\rightarrow A$ such that
\[
\left\| [\phi^{(i)}(x), y]\right\|\leq\varepsilon \left\| x\right\|
\]
for $i\in\{0, 1\}$, $x\in M_k$, $y\in\mathcal{F}$ and such that for each $\tau\in Y\cup Z$, there exists $i(\tau)\in\{0, 1\}$ such that $\tau(\phi^{(i(\tau))}(1_k))>1-\varepsilon$.
By Lemma \ref{4b}, we could reduce the number of cpc order zero maps to one, which is the $\phi$ we desired. 
\end{proof}

The recursive subhomogeneous algebras (RSH algebras) is an important class of C*-algebras. Recall the definition of a RSH algebra given in \cite{22}:

\begin{dfn}
A recursive subhomogeneous algebra is a C*-algebra given by the following recursive definition:\\
(1) If $X$ is a compact Hausdorff space and $n\geq 1$, then $C(X,M_n)$ is a recursive subhomogeneous algebra;\\
(2) If $A$ is a recursive subhomogeneous algebra, $X$ is a compact Hausdorff space, $X^{(0)}\subset X$ is closed, $\varphi: A\rightarrow C(X^{(0)},M_n)$ is any unital homomorphism, and $\rho: C(X,M_n)\rightarrow C(X^{(0)},M_n)$ is the restriction homomorphism, then the pullback
\[
A\oplus_{C(X^{(0)},M_n)} C(X,M_n)=\{(a,f)\in A\oplus C(X,M_n): \varphi(a)=\rho(f)\}.
\]
is a recursive subhomogeneous algebra.
\end{dfn}

\begin{prop}\label{4x}
Let $A$ be a simple nuclear separable unital infinite-dimensional C*-algebra with non-empty tracial state space. Suppose $T(A)$ is affinely homeomorphic to the tracial state space of some RSH algebra of finite topological dimension. The following conditions are equivalent:\\
(1) $A$ is $\mathcal{Z}$-stable;\\
(2) $A$ has strict comparison.
\end{prop}

\begin{proof}
We only need to construct a uniformly tracially large cpc order zero map.

First consider a zero-step RSH algebra of finite topological dimension, $C(X,M_n)$. This theorem follows immediately from Theorem 4.6 of \cite{14} because the extreme boundary of its tracial state space $X$ is compact and of finite topological dimension.

Then we consider a one-step RSH algebra of finite topological dimension:
\begin{align*}
B&=C(X,M_n)\oplus_{C(Y^{(0)},M_m)} C(Y,M_m)\\
&=\{(f,g)\in C(X,M_n)\oplus C(Y,M_m): \varphi(f)=\rho(g)\}
\end{align*}
where $X, Y$ are compact Hausdorff spaces of finite covering dimension, $Y^{(0)}\subset Y$ is closed, and $\varphi: C(X,M_n)\rightarrow C(Y^{(0)},M_m)$ is some unital homomorphism and $\rho: C(Y,M_m)\rightarrow C(Y^{(0)},M_m)$ is the restriction homomorphism. Since $T(B)$ and $T(A)$ are affinely homeomorphic, we claim that $\partial_eT(A)$ is of finite covering dimension and has the tightness property. Then the result follows by Theorem \ref{4a}. Indeed, $\partial_eT(B)$ is homeomorphic to $X\cup (Y\backslash Y^{(0)})$ whose cluster points in $X\cup Y$ all belongs to $Y^{(0)}$. But for any point $y$ in $Y^{(0)}$, which is the extreme point of $T(C(Y,M_m))$, $y\circ\varphi$ is a tracial state of $C(X,M_n)$. So there exists a Borel probability measure $\mu_y$ on $X$, such that $y\circ\varphi=\int_X d\mu_y$. Since $X$ is compact, this shows that $\partial_e T(B)$ has the tightness property, so does $\partial_e T(A)$.

Consider the case of a two step RSH algebra of finite topological dimension:
\begin{align*}
C&=B\oplus_{C(Z^{(0)},M_p)} C(Z,M_p)\\
&=\{(b,g)\in B\oplus C(Y,M_p): \varphi '(b)=\rho '(g)\}
\end{align*}
where $B$ is defined as above, $Z$ is a compact Hausdorff space of finite covering dimension, $Z^{(0)}\subset Z$ is closed, and $\varphi ': B\rightarrow C(Z^{(0)},M_p)$ is some unital homomorphism and $\rho ': C(Z,M_p)\rightarrow C(Z^{(0)},M_p)$ is the restriction homomorphism. Since $T(C)$ and $T(A)$ are affinely homeomorphic, and $\partial_eT(C)$ is homeomorphic to \mbox{$X\cup (Y\backslash Y^{(0)})\cup (Z\backslash Z^{(0)})$} whose cluster points belongs to $Y^{(0)}\cup Z^{(0)}$, so is $\partial_eT(A)$. By similar argument as the one-step case, for any point $y$ in $Y^{(0)}$, there exists a Borel probability measure $\mu_y$ on $X$, such that $y\circ\varphi=\int_X d\mu_y$. And for any point $z$ in $Z^{(0)}$, there exists a Borel probability measure $\mu_z$ on $X\cup (Y\backslash Y^{(0)})$, such that $z\circ\varphi=\int_{X\cup (Y\backslash Y^{(0)})} d\mu_z$. Hence by Definition 6.6, $Y\backslash Y^{(0)}$ is weakly tight relative to $X$, and $Z\backslash Z^{(0)}$ is weakly tight relative to $X\cup (Y\backslash Y^{(0)})$.

Since $X$ is compact and of finite covering dimension, by Lemma \ref{4d} and \ref{4y}, for each $k\geq 2$, there exists a cpc order map \mbox{$M_k\rightarrow A_\infty\cap A'$} which is uniformly tracially large on $X$. Then by Lemma \ref{4z}, there exists a cpc order map $M_k\rightarrow A_\infty\cap A'$ which is uniformly tracially large on $X\cup (Y\backslash Y^{(0)})$. Applying Lemma \ref{4z} again, we obtain a cpc order map $M_k\rightarrow A_\infty\cap A'$ which is uniformly tracially large on \mbox{$X\cup (Y\backslash Y^{(0)})\cup (Z\backslash Z^{(0)})$}.

Using the same type of arguments as in the two-step case, it is easy to see that for a general $n$-step RSH algebra of finite topological dimension, the extreme boundary $X$ of its tracial state space can be written as a union of subsets of $X$, say $X_1, X_2, ..., X_n$, such that the following holds: $X_1$ is compact; $X_2$ is weakly tight relative to $X_1$; $X_3$ is weakly tight relative to $X_1\cup X_2$; ...; $X_n$ is weakly tight relative to $\bigcup_{i=1}^{n-1}X_i$. Then for each $k\geq 2$, we first construct a cpc order zero map $M_k\rightarrow A_\infty\cap A'$ which is uniformly tracially large on $X_1$ using Lemma \ref{4d} and \ref{4y}, then by applying Lemma \ref{4z} $n-1$ times, finally we will obtain a cpc order zero map $M_k\rightarrow A_\infty\cap A'$ which is uniformly tracially large on $X$.
\end{proof}

\section*{Acknowledgment}
I would like to express my sincere appreciation to my advisor, Dr. Andrew Toms, for his support and guidance throughout the research. I would also like to thank the reviewers for their time and expertise, and for their valuable comments on the earlier version of this manuscript.

\end{document}